\documentclass[journal,twoside,web]{ieeecolor}
\usepackage{lcsys}
\usepackage{cite}
\usepackage{amsmath,amssymb,amsfonts}
\usepackage{mathtools}
\usepackage{algorithmic}
\usepackage{graphicx}
\usepackage{textcomp}
\usepackage{float}

\newcommand{\R}{\mathbb{R}}
\newcommand{\dd}{\,\mathrm{d}}

\newtheorem{theorem}{Theorem}[section]

\newtheorem{lemma}[theorem]{Lemma}
\newtheorem{corollary}[theorem]{Corollary}

\newtheorem{remark}{Remark}[section]

\newtheorem{assumption}{Assumption}[section]

\title{A Posteriori Second-Order Guarantees for Bolza Problems via Collocation}

\author{Dongzhe Zheng and Wenjie Mei$^{*}$%
\thanks{D. Zheng was with the Department of Computer Science and Engineering, School of Electronic Information and Electrical Engineering, Shanghai Jiao Tong University, Shanghai 200240, China, during the development of this work
        (e-mail: dz1011@wildcats.unh.edu, dz5992@princeton.edu).}%
\thanks{W. Mei is with the School of Robotics and Automation, Suzhou Campus,
        Nanjing University, Suzhou 215163, China
        (e-mail: mei.wenjie@nju.edu.cn).}%
\thanks{$^{*}$Corresponding author. This work was supported in part by the National Natural Science Foundation of China (NSFC)
under Grant 62403125, in part by the Natural Science Foundation of Jiangsu Province under Grant
BK20241283.}
}

\begin{document}

\maketitle
\thispagestyle{empty}
\pagestyle{empty}

%%%%%%%%%%%%%%%%%%%%%%%%%%%%%%%%%%%%%%%%%%%%%%%%%%%%%%%%%%%%%%%%%%%%%%%%%%%%%%%%
\begin{abstract}
Direct collocation for Bolza optimal control yields discrete Karush–Kuhn–Tucker (KKT) points, while practical solvers expose only discrete quantities such as primal-dual iterates, reduced Hessians, and Jacobians. This creates a gap between continuous second-order optimality theory and what can be certified from solver output. We develop an a posteriori certification framework that bridges this gap. Starting from a discrete KKT solution, we reconstruct piecewise polynomial state, control, and costate trajectories, evaluate residuals of the dynamics, boundary, and stationarity conditions, and derive a computable lower bound for the continuous second variation. The bound is expressed as the discrete reduced curvature minus explicit residual-dependent correction terms. A positive bound yields a sufficient certificate for continuous second-order sufficiency and provides quantitative information relevant to local growth and trust-region sizing. The constants entering the certification inequality are conservatively estimable from reconstructed discrete data. The resulting test is operationally verifiable from collocation outputs and naturally supports adaptive mesh refinement through residual decomposition. We also outline an extension to path inequalities with isolated transversal switches.
\end{abstract}


\section{Introduction}

\IEEEPARstart{C}{ontinuous-time} second-order sufficient conditions (SSOC) for Bolza optimal control are well understood: strong regularity, coercivity of the second variation, and quadratic growth have been thoroughly characterized in the variational literature \cite{Vinter2000,BonnansShapiro2000}. However, practical solvers based on direct collocation return only discrete nonlinear programming (NLP) data, as is typical in modern Gaussian-quadrature direct collocation transcriptions \cite{PattersonRao2014}, iterates $(x^N,u^N,\lambda^N)$, reduced Hessians, and KKT Jacobians, not continuous trajectories or functional second variations. 

This paper addresses the resulting gap by developing an a posteriori certification framework that constructs continuous optimality certificates from discrete solver output.

On the continuous side, we work in the KKT framework for the endpoint Lagrangian and invoke standard SSOC results for strong local optimality \cite{Vinter2000}. On the discrete side, we use the reduced Hessian of the collocation transcription as a computable curvature surrogate \cite{NocedalWright2006}, together with standard polynomial and costate reconstruction techniques \cite{Betts2010,FrancolinBensonHagerRao2015}. Recent work on Gauss and Radau collocation has focused primarily on convergence and stability of discrete solutions to continuous ones \cite{HagerHouRao2016,HagerLiuMohapatraRaoWang2018,HagerHouMohapatraRaoWang2019}; in contrast, we provide a posteriori certificates for continuous SSOC starting from an already-computed discrete KKT point.

Starting from a discrete KKT point, we reconstruct piecewise polynomial trajectories and evaluate residuals of dynamics, boundary, and stationarity. Combining a right-inverse bound for the linearized KKT mapping with smoothness bounds of the Lagrangian second derivatives yields a lower bound for the continuous second variation. The acceptance test reads
\begin{equation}\label{eq:acceptance-intro}
\hat{\alpha}_N\bigl(1-C_T E_N^{(2)}\bigr)^2 \;>\; \bigl(\Gamma + C_{\mathrm{quad}} + C_T'\bigr)\,E_N^{(2)},
\end{equation}
where $\hat{\alpha}_N$ is the minimum eigenvalue of the discrete reduced Hessian, $E_N^{(2)}$ collects $L^2$-residual measures, and the constants $C_T$, $\Gamma$, $C_{\mathrm{quad}}$, $C_T'$ are conservatively estimable from reconstructed discrete data. Positivity of the left-hand side minus the right-hand side provides a sufficient certificate for SSOC at the continuous level, together with quantitative bounds relevant to quadratic growth and trust-region sizing.

Algorithmically, the framework delivers an explicit acceptance test using only quantities available from the discrete solve: if the discrete curvature exceeds the residual-weighted threshold, the continuous problem is certified; otherwise, the mesh is refined or the polynomial degree is increased. The residual decomposition naturally guides adaptive refinement by identifying dominant error contributors \cite{PattersonHagerRao2015,LiuHagerRao2015,LiuHagerRao2018,MillerHagerRao2021}.

We also outline an extension to path inequalities with isolated transversal switches, where modified geometric constants preserve the same acceptance structure; see Section~\ref{sec:discussion} for discussion, with details deferred to future work.

\section{Problem Statement}\label{sec:problem}

We consider the Bolza problem with control $u\in L^2([0,T];\mathbb{R}^m)$ and no box or path inequality constraints in this section. The stationary point satisfies the interior optimality condition $H_u(t,x,u,p)=0$. The costate $p$ denotes the Lagrange multiplier for the state dynamics $\dot{x}=f(t,x,u)$. The second-order analysis is conducted for the endpoint Lagrangian via its second variation and SSOC, following the variational framework in \cite{Vinter2000}.

Let $T>0$. Consider $x\in W^{1,2}([0,T];\R^n)$, $u\in L^2([0,T];\R^m)$, we analyze the following optimal control problem:  
\begin{align}
\min_{x,u}\ J(x,u) &= K\big(x(0),x(T)\big) + \int_0^T L\big(t,x(t),u(t)\big)\,\dd t,  \notag \\
\dot x(t) &= f\big(t,x(t),u(t)\big), \quad t \in [0,T],  \label{eq:dyn}
\end{align}
where the functions $f,L,K\in C^2$ (see \cite{Vinter2000,BonnansShapiro2000} for the standard regularity assumptions in variational analysis).

\textbf{Problem Formulation:} The KKT condition defined below corresponds to the problem stated in \eqref{eq:dyn} with optional boundary equality constraints $b(x(0),x(T))=0$. Let $n_b\ge 0$ denote the number of boundary constraints and let $\lambda\in\mathbb{R}^{n_b}$ be the associated multiplier. When $n_b=0$ the terms involving $b$ and $\lambda$ vanish.

Define the Hamiltonian $H$ and the endpoint Lagrangian $\mathcal L$:
\begin{align}
H(t,x,\,u,\,p)&=L(t,x,u)+p^\top f(t,x,u), \label{eq:H_def} \\
\mathcal L(x,u,p,\lambda)&=K(x(0),x(T))+\lambda^\top b\big(x(0),x(T)\big)\notag\\
&\quad+\int_0^T \big(L(t,x,u)+p^\top(f(t,x,u)-\dot x)\big)\,\mathrm{d}t. \label{eq:L_def}
\end{align}
Define the KKT mapping $F=(F_1,\dots,F_6)$ by
\begin{align}
F_1&=\dot x-f(t,x,u)\;\in L^2, \label{eq:F1_def}\\
F_2&=-\dot p-H_x(t,x,u,p)\in L^2, \label{eq:F2_def}\\
F_3&=H_u(t,x,u,p)\in L^2, \label{eq:F3_def}\\
F_4&=b\big(x(0),x(T)\big)\in\mathbb{R}^{n_b}, \label{eq:F4_def}\\
F_5&=p(T)-K_{x_T}(x_0,x_T)-\big(b_{x_T}(x_0,x_T)\big)^\top\lambda\in\mathbb{R}^n, \label{eq:F5_def}\\
F_6&=-p(0)-K_{x_0}(x_0,x_T)-\big(b_{x_0}(x_0,x_T)\big)^\top\lambda. \label{eq:F6_def}
\end{align}
We equip the KKT mapping with the mixed norm
\begin{equation}\label{eq:KKT_norm}
\|F\|:=\|F_1\|_{L^2}+\|F_2\|_{L^2}+\|F_3\|_{L^2}+|F_4|+|F_5|+|F_6|,
\end{equation}
where $|\cdot|$ denotes the Euclidean norm on finite-dimensional vectors.
In the reconstruction, we additionally report an $L^\infty$-based diagnostic for the dynamics defect (see Section~\ref{sec:collocation}).

We let $T_{\mathrm{cont}}$ stand for the feasible tangent space of $(x,u)$-variations that satisfy the linearized dynamics and boundary conditions at a KKT point, and $T_{\mathrm{disc}}$ represents its discrete counterpart defined by the collocation scheme.

For $v=(\delta x,\delta u)\in T_{\mathrm{cont}}$, let $Q(x,u)(v)$ denote the second G\^ateaux derivative of the endpoint Lagrangian restricted to the feasible tangent space $T_{\mathrm{cont}}$ at $(x,u)$. 

If there exists $\alpha>0$ such that $Q(x^\star,u^\star)(v)\ge \alpha\|v\|^2$ for all $v\in T_{\mathrm{cont}}$, then $(x^\star,u^\star)$ is a strong local minimizer and exhibits quadratic growth in a neighborhood, a standard result for optimization problems whose defining functions are of class $C^2$, \emph{i.e.}, twice continuously differentiable \cite{BonnansShapiro2000}.

\begin{assumption}[Strengthened Legendre Condition]\label{ass:SLC} There exists $\rho>0$ such that along interior arcs (where no control bound is active),
\begin{equation}
H_{uu}(t,x,u,p)\succeq \rho I_{m \times m}
\end{equation}
holds almost everywhere.
\end{assumption}

Assumption~\ref{ass:SLC} ensures coercivity in the control direction and is required for the product-norm estimates in Section~\ref{sec:unconstrained}. For problems with active control bounds or singular arcs, the condition is replaced by second-order sufficiency on the critical cone; see \cite{OsmolovskiiMaurer2012}.

\section{Collocation, reconstruction, and residuals}\label{sec:collocation}

This section defines the reconstruction used in the a posteriori analysis, the associated residuals, and the stability constants that relate the continuous and discrete second-order objects. Throughout, constants such as $C_{\mathrm{geo}}$, $C_T$, $C_{\mathrm{quad}}$, $C_T'$, and $\Lambda$ denote computable upper bounds obtained from the discrete solution and standard interpolation estimates. These bounds are generally conservative but sufficient for the certification test.

\textbf{Reconstruction and residuals}: Let a Gauss, Radau, or Lobatto collocation scheme of degree $p$ on a mesh of size $h$ produce a discrete KKT point $(x^N,u^N)$ (and the associated discrete multipliers). From this discrete solution, we reconstruct piecewise polynomial trajectories $(X^N,u^N,p^N)$ for state, control, and costate that interpolate the discrete variables and satisfy the discrete adjoint relations, with a terminal condition
\begin{equation}
  p^N(T)
  = K_{x_T}\bigl(X^N(0),X^N(T)\bigr).
\end{equation}
The reconstruction $(X^N,u^N,p^N)$ is used to evaluate the continuous KKT equations.
Then, we distinguish two residual indicators. For theoretical perturbation and projection stability, define the $L^2$-aggregate
\begin{align}\label{eq:EN2_def}
&E_N^{(2)}
:= \bigl\|\dot{X}^N - f(\cdot,X^N,u^N)\bigr\|_{L^2} \notag\\
& + \bigl\|H_u(\cdot,X^N,u^N,p^N)\bigr\|_{L^2}
      + \bigl|b(X^N(0),X^N(T))\bigr|.
\end{align}
For diagnostics and visualization, define the $L^\infty$-indicator
\begin{align}\label{eq:Einf_def}
&E_\infty
:= \bigl\|\dot{X}^N - f(\cdot,X^N,u^N)\bigr\|_{L^\infty} \notag\\
& + \bigl\|H_u(\cdot,X^N,u^N,p^N)\bigr\|_{L^\infty}
      + \bigl|b(X^N(0),X^N(T))\bigr|.
\end{align}
These two are related by $\|r\|_{L^2}\le \sqrt{T}\|r\|_{L^\infty}$ for a function $r$, hence
$
E_N^{(2)}
\le \sqrt{T}\bigl(E_\infty-\bigl|b(X^N(0),X^N(T))\bigr|\bigr)
   + \bigl|b(X^N(0),X^N(T))\bigr|
\le \sqrt{T}\,E_\infty + \bigl|b(X^N(0),X^N(T))\bigr|$.
All theoretical bounds use $E_N^{(2)}$, while $E_\infty$ is used only as a diagnostic indicator.

Variations are equipped with the product norm
\begin{align}
  \|v\|^2
  &:= \|\delta x\|_{L^2}^2 + \|\delta u\|_{L^2}^2 + |\delta x(0)|^2 + |\delta x(T)|^2,\label{eq:var_norm}
\end{align}
where $\|\cdot\|$ denotes the product norm on the variation space
$L^2(0,T;\mathbb{R}^n)\times L^2(0,T;\mathbb{R}^m)\times \mathbb{R}^n\times \mathbb{R}^n$,
with $\|\delta x\|_{L^2}$, $\|\delta u\|_{L^2}$ the usual $L^2$ norms and
$|\delta x(0)|$, $|\delta x(T)|$ the Euclidean norms of the endpoint variations;
the endpoint terms are included to control boundary contributions to the second variation.

\textbf{Regularity and linearized KKT map}: We assume that $f$, $L$, and $K$ belong to $C^{2,1}$ on a neighborhood of $(X^N,u^N)$, \emph{i.e.}, they are twice continuously differentiable and their second derivatives are locally Lipschitz continuous. In this neighborhood, we introduce the uniform curvature bound {\small
\begin{equation}\label{eq:L2_def} 
  L_2 :=
  \sup \bigl\{
     \|\nabla^2 L(t,x,u)\|,
     \ \|\nabla^2 f_i(t,x,u)\|,
     \ \|\nabla^2 K(x_0,x_T)\|
  \bigr\},
\end{equation} }
where the supremum is taken over all admissible $(t,x,u,x_0,x_T)$ and all components $f_i$ of $f$.
The constant $L_2$ bounds the size of all second derivatives entering the endpoint Lagrangian and hence the magnitude of the continuous second variation, $\|\cdot\|$ denotes matrix induced norm.

Along the reconstruction, we set
\begin{align}
  A(t) &:= f_x\bigl(t,X^N(t),u^N(t)\bigr), \label{eq:Adef}\\
  B(t) &:= f_u\bigl(t,X^N(t),u^N(t)\bigr), \label{eq:Bdef}
\end{align}
so that $A(t)$ and $B(t)$ describe the linearized dynamics. Let $F$ denote the KKT mapping defined in~\eqref{eq:F1_def}--\eqref{eq:F6_def}, and let
\begin{equation}
  \mathcal{F}_N := DF(X^N,u^N,p^N,\lambda^N)
\end{equation}
be its Fr\'echet derivative at the reconstructed primal--dual point (with boundary multipliers $\lambda^N$ when present).

We assume that the KKT generalized equation is strongly regular, as demonstrated in~\cite{Robinson1980} (see also~\cite{DontchevRockafellar2009,BonnansShapiro2000}). Strong regularity implies that $\mathcal{F}_N$ admits a bounded right inverse, and we define the geometric constant
\begin{equation}\label{eq:Cgeo}
  C_{\mathrm{geo}} := \bigl\|\mathcal{F}_N^{-1}\bigr\|,
\end{equation}
where $\|\cdot\|$ denotes bounded linear operator norm.

\begin{remark}[Computable upper bound for $C_{\mathrm{geo}}$]\label{rem:Cgeo_compute}
Let $M_h$ denote the discrete KKT Jacobian obtained by linearizing the collocation equations at $(x^N,u^N,p^N,\lambda^N)$. Let $L_h$ and $R_h$ be the lifting and restriction operators between the continuous and discrete spaces. Then,
\begin{equation}\label{eq:Cgeo_bound}
C_{\mathrm{geo}} \;\le\; \|L_h\|\,\|M_h^{-1}\|_2\,\|R_h\| + \varepsilon_h,
\end{equation}
where spectral norm $\|M_h^{-1}\|_2$ is evaluated via singular-value decomposition and $\varepsilon_h = O(h^p)$ accounts for nonconformity (for trapezoidal collocation with piecewise linear reconstruction, $\|L_h\|\le 2$ and $\|R_h\|\le 1$).
\end{remark}

Combining $C_{\mathrm{geo}}$ with the curvature bound $L_2$ yields
\begin{equation}\label{eq:Gamma_def}
  \Gamma := C_{\mathrm{geo}} L_2,
\end{equation}
which weights the contribution of $E_N^{(2)}$ in the curvature transfer inequality.

To connect the computable discrete error quantity \(E_N^{(2)}\) with the distance to a nearby exact KKT point, we next present a standard perturbation result: Under local strong regularity of the KKT generalized equation, this yields the following estimate. Standard perturbation results for strongly regular generalized equations imply that any exact KKT point $(x^\star,u^\star,p^\star,\lambda^\star)$ close to $(X^N,u^N,p^N,\lambda^N)$ satisfies
\begin{equation}\label{eq:distance-EN}
  \|x^\star - X^N\|_{0,\infty} + \|u^\star - u^N\|_{L^2}
  \;\le\; C_\ast E_N^{(2)}
\end{equation}
for some constant $C_\ast>0$ depending only on local problem data, where $\|\cdot\|_{0,\infty}$ denotes the standard supremum norm.

\textbf{Tangent spaces, projection, and quadrature}: We first recall the definitions of $T_{\mathrm{cont}}$ and $T_{\mathrm{disc}}$, presented after Equation~\eqref{eq:KKT_norm}. 

The collocation basis defines a projection $\Pi_N : T_{\mathrm{cont}} \to T_{\mathrm{disc}}.$
By standard stability estimates for linear variational equations and stable collocation schemes (refer to, \emph{e.g.},~\cite{Betts2010,DontchevRockafellar2009,ChenDuHagerYang2019}), there exists a constant $C_T \ge 0$, independent of $h$, for all sufficiently small $h$, such that for all $v\in T_{\mathrm{cont}}$,
\begin{equation}\label{eq:near-isometry}
  (1 - C_T E_N^{(2)})\,\|v\|
  \;\le\; \|\Pi_N v\|
  \;\le\; (1 + C_T E_N^{(2)})\,\|v\|.
\end{equation}

\begin{lemma}[Projection stability constant]\label{lem:CT_bound}
If strong regularity and the condition~\eqref{eq:distance-EN} hold, then the linearized variational coefficients are perturbed by $O(E_N^{(2)})$. Applying Gronwall's inequality to the variational equation and bounding the interpolation operator yields a computable upper bound
\begin{equation}\label{eq:CT_bound}
C_T \;\le\; c_\Pi\,e^{\|A\|_{L^\infty} T}\Bigl(1 + \frac{\|B\|_{L^\infty}}{\rho} \Bigr),
\end{equation}
where $c_\Pi$ is the Lebesgue constant of the interpolation scheme (for example, for piecewise linear or cubic Hermite, $c_\Pi\le 2$), and $\rho$ is defined as in Assumption~\ref{ass:SLC}.
\end{lemma}

Let $Q(x,u)$ denote the continuous second variation of the endpoint Lagrangian restricted to $T_{\mathrm{cont}}$ (see Section~\ref{sec:problem}), and denote
\begin{equation}
  Q_N(v) := Q(X^N,u^N)(v), \qquad v\in T_{\mathrm{cont}}.
\end{equation}
Let $Q_{\mathrm{disc}}$ be the reduced quadratic form on $T_{\mathrm{disc}}$ induced by the collocation transcription (the discrete reduced Hessian). 

The projection, quadrature, and linearization errors can be collected in two additional constants:

i) The constant $C_{\mathrm{quad}}\ge0$ bounds the discrepancy between the collocation quadrature and the exact integral in $Q_N$. Let $v'$ and $v''$ denote the first and second derivatives of $v$ with respect to time, and $g_k := g(t_k)$, $g_{k+1} := g(t_{k+1})$. For trapezoidal quadrature on each subinterval $[t_k,t_{k+1}]$, the local error satisfies $|\int_{t_k}^{t_{k+1}}g \,dt - \tfrac{h}{2}(g_k+g_{k+1})|\le \tfrac{h^3}{12}\|g''\|_{L^\infty}$. Bounding $\|v'\|$ and $\|v''\|$ via the variational equation and $L_{2,1}$ (the Lipschitz constant of second derivatives) yields the estimate
\begin{equation}\label{eq:Cquad_bound}
C_{\mathrm{quad}} \;\le\; c_{\mathrm{trap}}\,h^2,
\end{equation}
where $c_{\mathrm{trap}} >0$ depends on $\|A\|_{L^\infty}$, $\|B\|_{L^\infty}$, $\rho$, and $L_{2,1}$.

ii) The constant $C_T'\ge0$ bounds the nonconformity errors due to projection and lifting between $T_{\mathrm{cont}}$ and $T_{\mathrm{disc}}$. For piecewise polynomial reconstruction of degree $p$, standard interpolation error estimates give $C_T' = O(h^p)$ with an explicit constant depending on $p$ and mesh regularity. \vspace{0.5em}

Then, there exist $C_{\mathrm{quad}},C_T'\ge0$, depending on the scheme and local bounds on $A,B$, such that for all $v\in T_{\mathrm{cont}}$, 
\begin{equation}\label{eq:Q-transfer}
  Q_N(v)
  \;\ge\;
  Q_{\mathrm{disc}}(\Pi_N v)
  - \bigl(\Gamma E_N^{(2)} + C_{\mathrm{quad}} E_N^{(2)} + C_T' E_N^{(2)}\bigr)\,\|v\|^2.
\end{equation} 
Inequalities~\eqref{eq:near-isometry} and~\eqref{eq:Q-transfer} are the key properties used in Section~\ref{sec:unconstrained} to transfer discrete curvature to the continuous second variation.

Finally, by the $C^{2,1}$-regularity and compactness of the neighborhood of interest, the map $(x,u)\mapsto Q(x,u)$ is locally Lipschitz. Hence, there exists $\Lambda>0$ such that for all $(x,u)$ satisfying
$\|(x-X^N,u-u^N)\|\le r$, where $\|\cdot\|$ is the product norm defined in~\eqref{eq:var_norm},
and all $v\in T_{\mathrm{cont}}$,
\begin{equation}\label{eq:Q-Lipschitz}
  \bigl|Q(x,u)(v) - Q_N(v)\bigr|
  \;\le\; \Lambda r \,\|v\|^2.
\end{equation}
This bound is combined with~\eqref{eq:Q-transfer} and the distance estimate~\eqref{eq:distance-EN} in Section~\ref{sec:unconstrained} to obtain an a posteriori lower bound on the continuous second variation at a KKT point.

\section{Unconstrained curvature transfer and a posteriori certification}\label{sec:unconstrained}

The reconstruction and residuals from Section~\ref{sec:collocation} quantify how well a discrete KKT point approximates a continuous solution. We now combine this information with the reduced Hessian of the collocation problem to obtain a computable lower bound on the continuous second variation and thus an a posteriori certificate of strong local optimality.

\subsection{Curvature transfer}

Let $Q_{\mathrm{disc}}$ be the discrete reduced quadratic form on $T_{\mathrm{disc}}$, associated with the reduced Hessian of the collocation nonlinear program (see \cite{NocedalWright2006}). Let $\hat\alpha_N$ be its smallest eigenvalue, and let $\Pi_N : T_{\mathrm{cont}}\to T_{\mathrm{disc}}$ be the projection induced by the collocation basis with stability constant $C_T>0$ from Section~\ref{sec:collocation}. Then, for all $v\in T_{\mathrm{cont}}$,
\begin{align}
Q_{\mathrm{disc}}(\Pi_N v)
&\;\ge\; \hat\alpha_N\,\|\Pi_N v\|^2 \notag\\
&\;\ge\; \hat\alpha_N\,\bigl(1-C_T E_N^{(2)}\bigr)^2\,\|v\|^2,
\label{eq:Qdisc-lb}
\end{align}
and the reconstruction, quadrature, and projection estimates yield~\eqref{eq:Q-transfer}. 
Here, $Q_N$ is the continuous second variation of the endpoint Lagrangian along $(X^N,u^N)$, restricted to $T_{\mathrm{cont}}$. The constants $\Gamma$, $C_{\mathrm{quad}}$, and $C_T'$ characterize the problem geometry, quadrature, and nonconformity errors are defined as in Section~\ref{sec:collocation}. 

Combining \eqref{eq:Qdisc-lb} and \eqref{eq:Q-transfer} gives, for all $v\in T_{\mathrm{cont}}$,
\begin{equation}\label{eq:Qlower_bound}
\begin{aligned}
Q_N(v)\;\ge\; 
& \Bigl(\hat\alpha_N\bigl(1-C_T E_N^{(2)}\bigr)^2 \\ 
&-\Gamma E_N^{(2)}
- C_{\mathrm{quad}}E_N^{(2)} - C_T' E_N^{(2)}\Bigr)\,\|v\|^2.
\end{aligned}
\end{equation}
We define the transferred curvature
\begin{equation}\label{eq:alphacont}
\alpha_{\mathrm{cont}}\;:=\;\hat\alpha_N\bigl(1-C_T E_N^{(2)}\bigr)^2-\Gamma E_N^{(2)} 
- C_{\mathrm{quad}}E_N^{(2)} - C_T' E_N^{(2)}.
\end{equation}
If $\alpha_{\mathrm{cont}}>0$, then \eqref{eq:Qlower_bound} becomes
\begin{equation}
Q_N(v)\;\ge\;\alpha_{\mathrm{cont}}\|v\|^2\qquad\forall\,v\in T_{\mathrm{cont}},
\end{equation}
in terms of computable discrete quantities.

\subsection{Lipschitz stability and quadratic growth}

We formalize the Lipschitz stability of the second variation and derive computable bounds for the trust-region radius.

\begin{lemma}[A posteriori proximity bound]\label{lem:proximity}
Define the $L^\infty$-residuals
\begin{align}
e_{\mathrm{dyn},\infty} &:= \|\dot{X}^N - f(\cdot,X^N,u^N)\|_{L^\infty}, \label{eq:edyn_inf}\\
e_{\mathrm{adj},\infty} &:= \|-\dot{p}^N - H_x(\cdot,X^N,u^N,p^N)\|_{L^\infty}, \label{eq:eadj_inf}\\
e_{\mathrm{stat},\infty} &:= \|H_u(\cdot,X^N,u^N,p^N)\|_{L^\infty}, \label{eq:estat_inf}
\end{align}
and let $e_{\mathrm{bc}}$ collect the boundary residuals from~\eqref{eq:F4_def}--\eqref{eq:F6_def}. Here we slightly strengthen the earlier definition of $E_\infty$ for the KKT system by including the adjoint residual, i.e., $E_\infty := e_{\mathrm{dyn},\infty} + e_{\mathrm{adj},\infty} + e_{\mathrm{stat},\infty} + e_{\mathrm{bc}}$. Under Assumptions~\ref{ass:SLC} and strong regularity, if $E_\infty$ is sufficiently small, there exists a unique exact KKT point $(x^\star,u^\star,p^\star,\lambda^\star)$ satisfying
\begin{align}\label{eq:proximity_bound}
&\|x^\star - X^N\|_{0,\infty} + \|u^\star - u^N\|_{0,\infty} + \|p^\star - p^N\|_{0,\infty} \notag \\ & \le\; C_{\mathrm{close},\infty}\,E_\infty,
\end{align}
where
\begin{equation}\label{eq:Cclose_def}
C_{\mathrm{close},\infty} := C_{xp,\infty} + C_{u,\infty},
\end{equation}
with $C_{xp,\infty} \lesssim C_{\mathrm{geo}}(1+T)\exp\bigl((\|A\|_{L^\infty}+\|B\|_{L^\infty})T\bigr)$ and $C_{u,\infty} := \rho^{-1}(\|H_{ux}\|_\infty + \|H_{up}\|_\infty)C_{xp,\infty} + \rho^{-1}$.
\end{lemma}

\begin{lemma}[Computable Second-Variation Lipschitz bound]\label{lem:Lip_Q}
Define the tubular neighborhood $\Omega := \{(t,x,u,p) : t\in[0,T],\, \|x-X^N(t)\|\le\Delta_x,\, \|u-u^N(t)\|\le\Delta_u,\, \|p-p^N(t)\|\le\Delta_p\}$. Let
\begin{align}
L_{2,1}^H &:= L_{2,1}^L + P_{\max} L_{2,1}^f, \label{eq:L21H}\\
\Lambda &:= C_{\mathrm{int}}\bigl(L_{2,1}^H + M_{2,f}\bigr) + 2L_{2,1}^K, \label{eq:Lambda_def}
\end{align}
where $L_{2,1}^L$, $L_{2,1}^f$, $L_{2,1}^K$ are the Lipschitz constants of the second derivatives of $L$, $f$, $K$ on $\Omega$, $M_{2,f} := \sup_\Omega \max_i \|\nabla^2 f_i\|$, $P_{\max} := \|p^N\|_{0,\infty} + \Delta_p$, and $C_{\mathrm{int}} := \max\{T,1\}$. Then, for all $(x,u,p)\in\Omega$ and $v\in T_{\mathrm{cont}}$,
\begin{equation}\label{eq:Lip_Q_bound}
\bigl|Q(x,u)(v) - Q_N(v)\bigr| \;\le\; \Lambda\,\mathcal{R}(z;z_N)\,\|v\|^2,
\end{equation}
where $\mathcal{R}(z;z_N) := \|x-X^N\|_{0,\infty} + \|u-u^N\|_{0,\infty} + \|p-p^N\|_{0,\infty}$.
\end{lemma}

If $\mathcal{R}(z;z_N)\le r$, where
\begin{equation}\label{eq:trust-radius}
r:=\frac{\alpha_{\mathrm{cont}}}{2\Lambda},
\end{equation}
then, by combining \eqref{eq:Qlower_bound} and \eqref{eq:Lip_Q_bound}, we obtain
\begin{align}
Q(x,u)(v)
&\ge Q_N(v)-\Lambda \mathcal{R}(z;z_N)\,\|v\|^2 \notag\\
&\ge Q_N(v)-\Lambda r\,\|v\|^2 \notag\\
&\ge \bigl(\alpha_{\mathrm{cont}}-\Lambda r\bigr)\|v\|^2
= \tfrac12 \alpha_{\mathrm{cont}}\|v\|^2,
\label{eq:strong_convexity}
\end{align}
for all $v\in T_{\mathrm{cont}}$.
Hence, if $\alpha_{\mathrm{cont}}>0$, then $Q(x,u)$ satisfies a uniform second-order sufficient condition on $T_{\mathrm{cont}}$ throughout the neighborhood
$
\{\, z:\mathcal{R}(z;z_N)\le r \,\},
$
which in turn implies quadratic growth; see, \emph{e.g.}, \cite{Vinter2000}.

\begin{corollary}[A posteriori certified quadratic growth]\label{cor:certification}
Define $\Gamma_{\mathrm{tot}} := \Gamma + C_{\mathrm{quad}} + C_T'$ and consider  $\alpha_{\mathrm{cont}}$ defined in~\eqref{eq:alphacont}. Accept the discrete solution if
\begin{equation}\label{eq:exact_acceptance}
\hat{\alpha}_N\bigl(1 - C_T E_N^{(2)}\bigr)^2 \;>\; \Gamma_{\mathrm{tot}}\,E_N^{(2)}.
\end{equation}
If additionally $C_T E_N^{(2)} \le 0.1$, one may use the simplified test $\hat{\alpha}_N > \Gamma_{\mathrm{tot}} E_N^{(2)}$. When~\eqref{eq:exact_acceptance} holds and Lemma~\ref{lem:proximity} verifies $C_{\mathrm{close},\infty} E_\infty \le r$, the exact KKT point $(x^\star,u^\star)$ satisfies
\begin{equation}
Q(x^\star,u^\star)(v) \;\ge\; \tfrac{1}{2}\alpha_{\mathrm{cont}}\|v\|^2, \quad \forall\,v\in T_{\mathrm{cont}}.
\end{equation}
\end{corollary}

\subsection{Practical a posteriori check and acceptance criterion}

In practice, $\hat\alpha_N$ is the smallest eigenvalue of the discrete reduced Hessian, and $E_N^{(2)}$ is the residual indicator from Section~\ref{sec:collocation}. For the chosen mesh and reconstruction, we estimate
$
\Gamma_{\mathrm{tot}}
$. 
The certification reduces to the scalar test~\eqref{eq:exact_acceptance}. If this holds, then $\alpha_{\mathrm{cont}}>0$; we accept the discrete solution, evaluate $\alpha_{\mathrm{cont}}$ via~\eqref{eq:alphacont}, and set the trust-region radius $r=\alpha_{\mathrm{cont}}/(2\Lambda)$ from~\eqref{eq:trust-radius}. Otherwise, the mesh or polynomial degree is refined, and the procedure is repeated.

\section{Numerical example: Planar quadrotor maneuver}\label{sec:examples}

We illustrate the a posteriori certification framework on a nonlinear planar quadrotor benchmark. All geometric and stability constants are estimated from the discrete data, demonstrating the constructive nature of the framework.

\subsection{Dynamics and their optimal control}\label{subsec:quadrotor-formulation}

Consider the planar rigid-body quadrotor with state $x = [y,\, z,\, \theta,\, v_y,\, v_z,\, \omega]^\top \in \mathbb{R}^6$ (horizontal and vertical position, pitch angle, and their rates) and control $u = [u_1,\, u_2]^\top \in \mathbb{R}^2$ (left and right rotor thrusts). The system parameters are mass $m = 1.0$ kg, gravitational acceleration $g = 9.81$ m/s$^2$, half-arm length $l = 0.3$ m, and moment of inertia $J_x = 0.2$ kg$\cdot$m$^2$.

The nonlinear dynamics are
\begin{align}
\dot{y} &= v_y, \quad \dot{z} = v_z, \quad \dot{\theta} = \omega, \quad \dot{v}_y = -\frac{u_1+u_2}{m}\sin\theta, \label{eq:quad_kin}\\
\dot{v}_z &= \frac{u_1+u_2}{m}\cos\theta - g, \quad \dot{\omega} = \frac{l}{J_x}(u_1 - u_2). \label{eq:quad_trans}
\end{align}

Over the horizon $[0,T]$ with $T = 2$s, we minimize
\begin{align}\label{eq:quad_cost}
J(x,u) &= \frac{1}{2}\int_0^T \bigl((x-x_{\mathrm{ref}})^\top Q(x-x_{\mathrm{ref}})  \\ & + u^\top R u\bigr)\,dt +
\frac{1}{2}(x(T)-x_f)^\top K(x(T)-x_f),
\end{align}
where $Q = \mathrm{diag}(1,1,0.1,0.1,0.1,0.1)$, $R = \mathrm{diag}(0.01,0.01)$, and the terminal penalty $K = 100 I_{6 \times 6}$. The reference $x_{\mathrm{ref}}(t) \equiv 0$, initial state $x_0 = [-1,\, 0,\, 0,\, 0,\, 0,\, 0]^\top$, and target $x_f = [1,\, 0.5,\, 0,\, 0,\, 0,\, 0]^\top$ define a displacement-to-hover maneuver.

Since the dynamics are affine in $u$, we have $H_{uu} = L_{uu} = R$. Thus Assumption~\ref{ass:SLC} holds globally with $\rho = \lambda_{\min}(R) = 0.01$.

\subsection{Discretization and reconstruction}

We discretize using Hermite--Simpson collocation on a uniform mesh with $N$ intervals. For each $N \in \{10, 15, 20, 25, 30, 35\}$, we solve the resulting NLP via SQP (tolerance $10^{-12}$) to obtain the discrete KKT point $(x^N, u^N, p^N)$, then reconstruct continuous trajectories $(X^N, u^N, p^N)$ via cubic Hermite interpolation for states and piecewise linear interpolation for controls.

The residuals are evaluated as in Section~\ref{sec:collocation}:
\begin{align}
e_{\mathrm{dyn},\infty} &= \|\dot{X}^N - f(\cdot, X^N, u^N)\|_{L^\infty}, \\
e_{\mathrm{stat},\infty} &= \|H_u(\cdot, X^N, u^N, p^N)\|_{L^\infty}, \\
e_{\mathrm{bc}} &= \|X^N(0) - x_0\| + \|X^N(T) - x_f\|_K,
\end{align}
with $E_\infty = e_{\mathrm{dyn},\infty} + e_{\mathrm{stat},\infty} + e_{\mathrm{bc}}$ and $E_N^{(2)}$ obtained by integrating the squared residuals.

\subsection{Estimated constants and certification}

For the mesh with $N = 35$, we estimate all constants from Section~\ref{sec:unconstrained} using automatic differentiation along the reconstructed trajectory. The tubular neighborhood radii are set to $\Delta_x = \Delta_u = \Delta_p = 0.1$.

The discrete KKT Jacobian $M_h$ is assembled at $(x^N, u^N, p^N)$. Computing its singular values yields $\sigma_{\min} = 1.87 \times 10^{-2}$, while evaluating the Jacobian mapping via automatic differentiation gives the upper bound $C_{\mathrm{geo}} \le 53.39$.

Sampling second derivatives via automatic differentiation over $\Omega$ yields $M_{2,f} = 17.90$ and $L_{2,1}^f = 1.23$. From the aggregated constant analysis, we obtain the bound $\Lambda \le 1.09$.

The evaluated residuals are
$
E_N^{(2)} = 3.27 \times 10^{-14}, \quad E_\infty = 7.05 \times 10^{-14}.
$
The discrete reduced Hessian has a minimum eigenvalue
$
\hat{\alpha}_N = 6.29 \times 10^{-4}.
$

With the upper bound $\Gamma \le 979.47$, the residual threshold is estimated as $\hat{\alpha}_N^{th} = 4.67 \times 10^{-11}$. Since $\hat{\alpha}_N \gg \hat{\alpha}_N^{th}$, the transferred continuous curvature is strictly positive:
$
\alpha_{\mathrm{cont}} = \hat{\alpha}_N - \hat{\alpha}_N^{th} \approx 6.29 \times 10^{-4} > 0.
$
The trust-region radius is $r = \frac{\alpha_{\mathrm{cont}}}{2\Lambda} = 2.88 \times 10^{-4}.$

Finally, by Lemma~\ref{lem:proximity}, with the bound $C_{\mathrm{close},\infty} \le 59.36$, we have
$
C_{\mathrm{close},\infty} E_\infty = 59.36 \times 7.05 \times 10^{-14} = 4.18 \times 10^{-12}.
$
Since $4.18 \times 10^{-12} \ll r$, the proximity condition $C_{\mathrm{close},\infty} E_\infty \le r$ is satisfied. Thus, certification is achieved even at relatively sparse discretizations.

\begin{figure}[h!]
  \centering
  \includegraphics[width=0.85\linewidth]{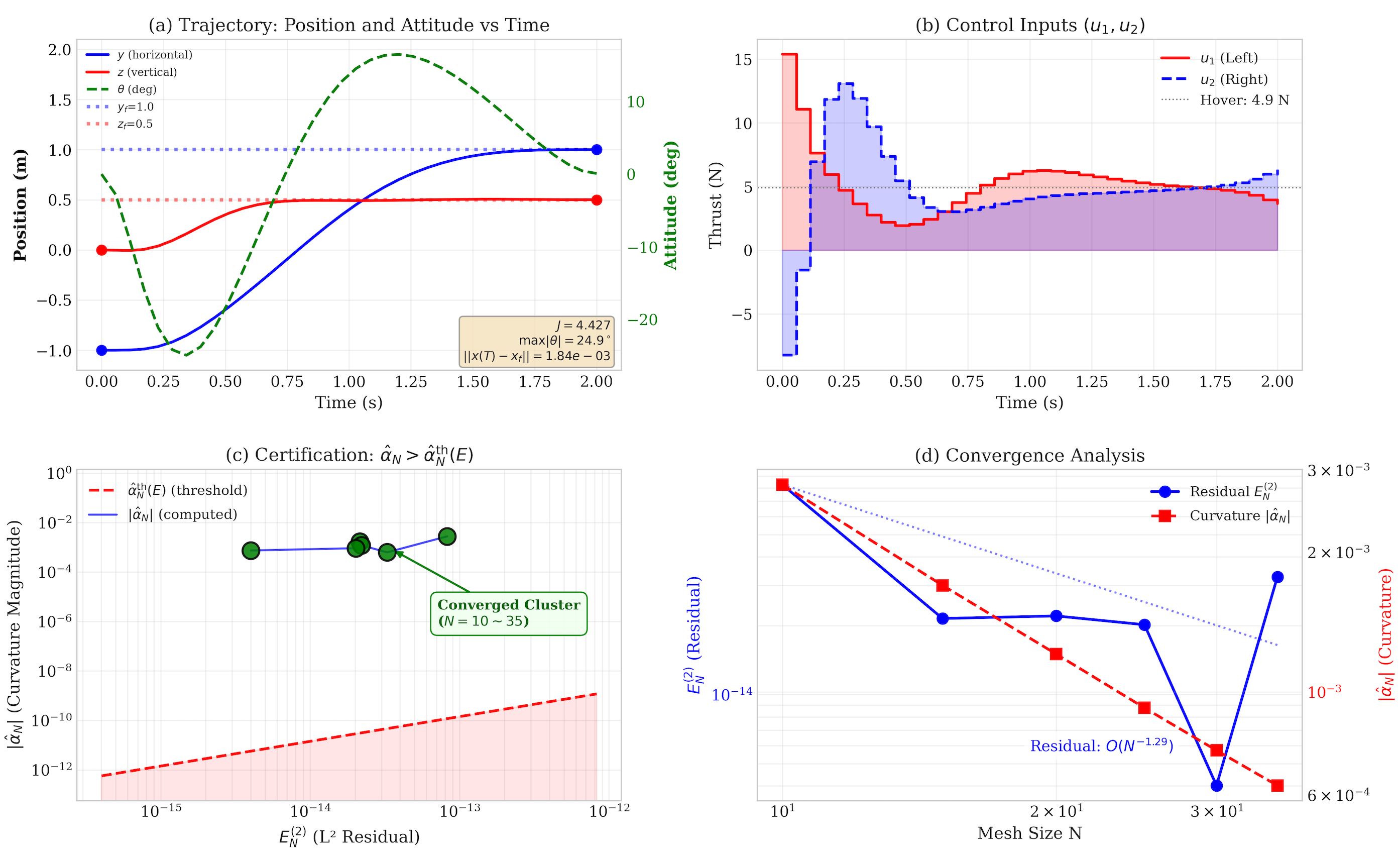}
  \caption{Planar quadrotor maneuver. (a) Reconstructed position $(y,z)$ and attitude $\theta$ on the certified mesh ($N=35$). (b) Rotor thrusts $u_1$ and $u_2$. (c) Discrete reduced curvature magnitude $|\hat{\alpha}_N|$ and residual threshold $\hat{\alpha}_N^{th}(E)$ versus the $L^2$-residual $E_N^{(2)}$. (d) Convergence of the $L^2$-residual and curvature magnitude with respect to the mesh size $N$. For all evaluated mesh sizes $N \ge 10$, the precise gradients provided by automatic differentiation drive the residuals near machine epsilon, satisfying the certification conditions $\hat{\alpha}_N > \hat{\alpha}_N^{th}$ uniformly.}
  \label{fig:quadrotor}
\end{figure}

\section{Discussion and Conclusion}\label{sec:discussion}

\subsection{An extension: isolated path-inequality switches}\label{subsec:extension}

Although the preceding framework and numerical validation address unconstrained controls, the structural advantage of our a posteriori approach is its extensibility to constrained topologies. We briefly discuss how the certification structure accommodates isolated path-inequality switches, deferring full derivations and numerical validation to future work. Recent direct methods with explicit switch and event detection point in a similar numerical direction \cite{NurkanovicSperlAlbrechtDiehl2024}.

Consider a path inequality $c(t,x,u)\le 0$ with $c:[0,T]\times\mathbb{R}^n\times\mathbb{R}^m\to\mathbb{R}$ of class $C^{2,1}$. When the active set consists of finitely many isolated, transversal crossings $\{\tau_k\}_{k=1}^K$ satisfying $\Delta\tau_{\min}:=\min_{k}(\tau_{k+1}-\tau_k)>0$, the path multiplier becomes a nonnegative atomic measure
\begin{equation}\label{eq:atomic_mu}
\mu = \sum_{k=1}^{K}\mu_k\,\delta_{\tau_k}, \qquad \mu_k \ge \mu_{\min} > 0.
\end{equation}
At each switching time $\tau_k$, the costate exhibits a jump:
\begin{equation}\label{eq:adjoint_jump}
p(\tau_k^+) - p(\tau_k^-) = \mu_k\,c_x\bigl(\tau_k,X(\tau_k),u(\tau_k)\bigr).
\end{equation}
In the discrete reconstruction of Section~\ref{sec:collocation}, these jumps are visible in $p^N$. Switching times are identified from the reconstructed slack $s^N(t):=c(t,X^N(t),u^N(t))$ by local search. The only switch-specific perturbations entering the stability constants are event-time errors and mass mismatches.

The reconstruction, residuals, and curvature-transfer pipeline established in Sections~\ref{sec:collocation}--\ref{sec:unconstrained} remain fundamentally useful. Rather than necessitating derivation from first principles, transversal switches act as structured perturbations to the linearized KKT map; referring to~\cite{OsmolovskiiMaurer2012,Vinter2000}, this effect can be encapsulated by inflating the geometric constant with switch-local diagnostics. Analytically, adjoint jumps manifest as low-rank block perturbations within the linearized KKT generalized equation. Under strong regularity and transversality~\cite{Robinson1980,DontchevRockafellar2009,OsmolovskiiMaurer2012}, standard sensitivity analysis indicates that the bounded right inverse experiences a controlled inflation relative to $C_{\mathrm{geo}}$. We isolate this topological effect into a modular scalar $\Gamma_{\mathrm{switch}}\ge 0$, yielding the updated curvature weight
\begin{equation}\label{eq:Gamma_tilde}
\widetilde{\Gamma} := \Gamma + \Gamma_{\mathrm{switch}}.
\end{equation}
The $L^2$-aggregate $E_N^{(2)}$ from Section~\ref{sec:collocation} is retained; any complementarity or slack diagnostic around switches is absorbed into $\Gamma_{\mathrm{switch}}$. The acceptance test becomes
\begin{equation}\label{eq:accept_switch}
\hat{\alpha}_N\bigl(1 - C_T E_N^{(2)}\bigr)^2 > \widetilde{\Gamma}\,E_N^{(2)},
\end{equation}
where $\hat{\alpha}_N$, $C_T$, and $E_N^{(2)}$ are defined as in Section~\ref{sec:unconstrained}. When~\eqref{eq:accept_switch} holds, the transferred continuous curvature remains positive and the trust-region construction applies. The ingredients to bound $\Gamma_{\mathrm{switch}}$ are embedded within $(X^N,u^N,p^N)$; its estimation follows from combining the perturbation calculus for strongly regular generalized equations~\cite{Robinson1980,DontchevRockafellar2009} with the second-order jump conditions for broken extremals~\cite{OsmolovskiiMaurer2012,Vinter2000}.

\subsection{Conclusion}\label{subsec:conclusion}

This paper established a constructive a posteriori certification framework linking discrete collocation solutions to continuous second-order sufficient conditions. Under strong regularity and the strengthened Legendre condition, all geometric and stability constants admit conservative but computable upper bounds derived from discrete data, yielding mesh refinement guidance and quantified trust-region guarantees. The bounds are generally conservative due to the use of global interpolation and Lipschitz estimates, but remain sufficient for the certification test.

\bibliographystyle{IEEEtran}
\bibliography{refs}

\end{document}